\documentclass[graybox]{svmult}
\usepackage{mathptmx}
\usepackage{helvet}
\usepackage{courier}
\usepackage{type1cm}
\usepackage{makeidx}
\usepackage{graphicx}
\usepackage{multicol}
\usepackage[bottom]{footmisc}
\usepackage{amssymb,amsmath,bm,bbm}
\begin{document}

\title*{Sparse Representations for Uncertainty Quantification of a Coupled Field-Circuit Problem}
\titlerunning{Uncertainty Quantification of a Coupled Field-Circuit Problem}
\author{Roland Pulch \and Sebastian Sch\"ops}
\institute{R.~Pulch: Universit\"at Greifswald, Institute of Mathematics and Computer Science, Walther-Rathenau-Str.~47, 17489 Greifswald, Germany, \email{roland.pulch@uni-greifswald.de}
\and
S.~Sch\"ops: Technische Universit\"at Darmstadt, Centre for Computational Engineering, Dolivostr.~15, 64293 Darmstadt, Germany, \email{schoeps@temf.tu-darmstadt.de}}

\maketitle
\abstract*{}
\abstract{We consider a model of an electric circuit,
  where differential algebraic equations for a circuit part are coupled
  to partial differential equations for an electromagnetic field part.
  An uncertainty quantification is performed by changing physical parameters
  into random variables.
  A random quantity of interest is expanded into the (generalised)
  polynomial chaos using orthogonal basis polynomials.
  We investigate the determination of sparse representations,
  where just a few basis polynomials are required for a sufficiently
  accurate approximation.
  Furthermore, we apply model order reduction with proper orthogonal
  decomposition to obtain a low-dimensional representation
  in an alternative basis.}


\section{Introduction}
\label{sec:introduction}
In science and engineering,
complex applications require an advanced modelling by multiphysics systems
or coupled systems.
We examine a coupled field-circuit problem of an electric circuit,
where differential algebraic equations (DAEs) for circuit components
are combined with partial differential equations (PDEs)
for electromagnetic components, see~\cite{pulch:schoeps-phd}.

Uncertainty quantification (UQ) investigates the impact of variations in
input parameters on a quantity of interest (QoI).
Often parameters are remodelled into random variables.
The random QoI can be expanded in the (generalised) polynomial chaos,
where orthogonal basis polynomials are involved, see~\cite{pulch:xiu-book}.
Sparse representations aim for a reduced set of basis functions
with a given accuracy of approximation.
Many methods for sparse representations have been derived and studied,
see~\cite{pulch:blatman,pulch:doostan,pulch:jakeman2017}
and the references therein. 
Alternatively, methods of model order reduction (MOR) yield low-dimensional 
(dense) approximations of the random QoI,
see~\cite{pulch:matcom2018,pulch:arxiv}.

We apply this UQ concept to the coupled field-circuit problem 
{\cite{pulch:trends}}.
On the one hand, sparse representations are determined by neglecting
basis functions with small coefficients.
On the other hand, MOR using proper orthogonal decomposition (POD) 
identifies a low-dimensional approximation in an alternative basis.
Our aim is to obtain approximations with as few basis functions
as possible, while still maintaining some accuracy.
Computational effort during the offline phase, i.e., evaluations of the
multiphysics systems, is not saved by the proposed methods. 
However, the online evaluation cost of the polynomials can be reduced
in the first approach.


\begin{figure}
  \begin{center}
    \includegraphics[width=8cm]{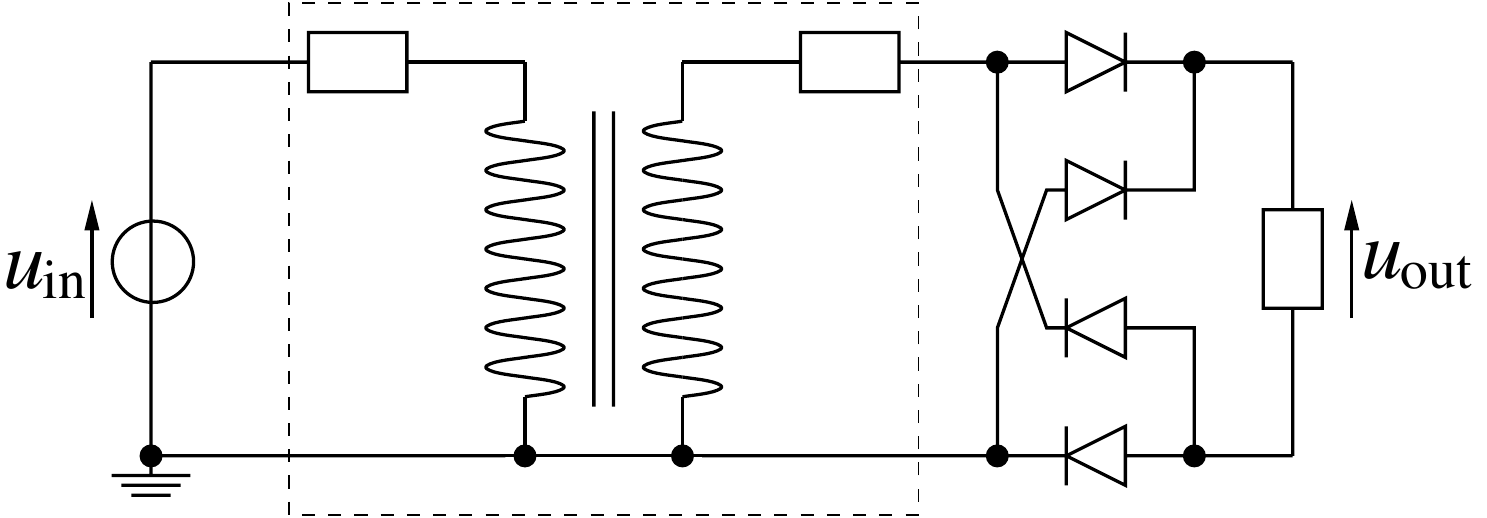}
  \end{center}
  \caption{Diagram of rectifier circuit.
    A PDE model is used for the components in dashed box.}
  \label{pulch:fig:circuit}
\end{figure}

\section{Coupled Field-Circuit Problem}
\label{sec:problem}
We investigate the rectifier circuit depicted in Fig.~\ref{pulch:fig:circuit}.
The model consists of a circuit part and a field part.
Modified nodal analysis (MNA) \cite{pulch:ho-etal} produces a system of DAEs 
\begin{equation} \label{pulch:mna}
  \begin{array}{rcl}
{\bf A}_C \frac{{\rm d}}{{\rm d}t} {\bf q}_C({\bf u},t) + 
{\bf A}_R {\bf r}({\bf u},t) + 
{\bf A}_L {\bf j}_L  
+ {\bf A}_V {\bf j}_V + 
{\bf A}_M {\bf j}_M  +
{\bf A}_D {\bf j}_D  +
{\bf A}_I {\bf i}(t)
& = & {\bf 0} , \\
\frac{{\rm d}}{{\rm d}t} \boldsymbol{\phi}_L({\bf j}_L,t) - 
{\bf A}_L^\top {\bf u} &  = & {\bf 0} , \\ 
{\bf A}_V^\top {\bf u} - {\bf v}(t) & = & {\bf 0} , \\
  \end{array}
\end{equation}
with incidence matrices ${\bf A}_\star$, node voltages ${\bf u}(t)$, branch currents ${\bf j}_L(t)$, ${\bf j}_V(t)$, sources ${\bf i}(t)$, ${\bf v}(t)$ and constitutive relations ${\bf q}_C(\cdot,t)$, ${\bf r}_C(\cdot,t)$, $\boldsymbol{\phi}_L(\cdot,t)$. Initial values are considered in the time interval $t \in [t_0,t_{\rm end}]$.
We apply Shockley's model
\begin{equation} \label{pulch:shockley}
  j_{D,k} = I_{S,k} \left( \exp \left( \textstyle {\bf A}_{D,k}^\top{\bf u} / U_{{\rm TH},k} \right) - 1 \right) , \qquad k=1,2,3,4
\end{equation}
for the four diodes with parameters $I_{S,k},U_{{\rm TH},k}$, where ${\bf A}_{D,k}$ denotes the $k$th column of ${\bf A}_{D}$. We involve a refined model for the transformer (dashed box in Fig.~\ref{pulch:fig:circuit}) given
by the two-dimensional (2D) magnetostatic approximation of Maxwell's equations
\begin{equation} \label{pulch:pde}
  \begin{array}{rcl}
  \nabla \cdot \left( 
\boldsymbol{\nu}(\| \nabla A(t,{\bf x}) \| ,{\bf x}) \; 
\nabla A (t,{\bf x}) \right) &=& \boldsymbol{\chi}({\bf x})^\top \mathbf{j}_M(t)
\qquad {\rm for}\;\; {\bf x} \in \Lambda \\
\displaystyle \frac{{\rm d}}{{\rm d}t}
\int_{\Lambda} \boldsymbol{\chi}({\bf x}) \; A(t,{\bf x}) \;
\mathrm{d}\mathbf{x}
&=&{\bf A}_M^\top {\bf u}(t) \\
\end{array}
\end{equation}
on the spatial domain~$\Lambda \subset \mathbbm{R}^2$.
The magnetic vector potential 
$A: [t_0,t_{\rm end}] \times \Lambda \rightarrow \mathbbm{R}$
is unknown.
The lumped currents and voltages are distributed and integrated
by the winding function 
$\boldsymbol{\chi} : \Lambda \rightarrow \mathbbm{R}^2$,
see~\cite{pulch:winding}.
The reluctivity is 
$\boldsymbol{\nu} :
\mathbbm{R} \times \Lambda \rightarrow \mathbbm{R}^{2 \times 2}$.
In the iron core region, it reads as
$\boldsymbol{\nu}(B,\vec{x}) = \nu(B) {\bf I}_2$ (identity matrix~${\bf I}_2$)
using Brauer's model 
\begin{equation} \label{pulch:brauer}
  \nu(B) = \kappa_1 \exp \left( \kappa_2 B^2 \right) + \kappa_3
\end{equation}
with the magnetic field~$B=\|\nabla A\|$ and the parameters
$\kappa_1,\kappa_2,\kappa_3$.
A finite element method yields a nonlinear system of algebraic equations.
More details on this coupled problem can be found in~\cite{pulch:schoeps-phd}.
Now we define the output voltage as QoI.


\section{Stochastic Model}
\label{sec:stochastic}
We consider uncertainties in $q=11$~parameters:
the parameters of Shockley's model~(\ref{pulch:shockley})
for each diode separately (8~parameters) and
the three parameters of Brauer's model~(\ref{pulch:brauer}).
We describe the uncertainties by independent uniform probability
distributions with 20\% variation around each mean value.
The random variables are ${\bf p}: \Omega \rightarrow \Pi$ with
event space~$\Omega$ and parameter domain $\Pi \subset \mathbbm{R}^{q}$.
The joint probability density function is constant on
the cuboid~$\Pi$.
Let $y : [t_0,t_{\rm end}] \times \Pi \rightarrow \mathbbm{R}$ be
the random output voltage (QoI) of the coupled problem.

The expected value of a function
$f : \Pi \rightarrow \mathbbm{R}$ reads as
\begin{equation} \label{pulch:expected-value}
  \mathbb{E} [f] = \displaystyle
  \frac{1}{{\rm volume}(\Pi)}
  \int_{\Pi} f({\bf p}) \; {\rm d}{\bf p} .
\end{equation}

The expected value~(\ref{pulch:expected-value}) implies an inner product
$< f , g > = \mathbb{E} [fg]$ for two square-integrable functions.
The accompanying norm is
$\| f \|_{{L}^2} = \sqrt{< f,f >}$.
We define the basis polynomials $( \Phi_i )_{i \in \mathbbm{N}}$
with $\Phi_i : \Pi \rightarrow \mathbbm{R}$ by
$ \Phi_i({\bf p}) =
\phi_{i_1} (p_1) \phi_{i_2} (p_2) \cdots \phi_{i_q} (p_q) $,
where $\phi_{\ell}$ denotes the (normalised) Legendre polynomial
of degree~$\ell$.
There is a one-to-one mapping between the indices~$i$ and
the multi-indices $i_1,\ldots,i_q$.
It follows that $( \Phi_i )_{i \in \mathbbm{N}}$ represents a complete
orthonormal system satisfying
$< \Phi_i , \Phi_j > = \delta_{ij}$.

We assume that the random process $y(t,\cdot)$ is square-integrable
for each~$t$.
Consequently, the (generalised) polynomial chaos expansion
\begin{equation} \label{pulch:gpc}
  y(t,{\bf p}) = \sum_{i=1}^\infty w_i(t) \Phi_i({\bf p})
\end{equation}
exists pointwise for each~$t$.
The coefficient functions are given by the inner products 
$w_i(t) = < y (t,\cdot) , \Phi_i(\cdot) >$.
The infinite series~(\ref{pulch:gpc}) is truncated to a finite sum
\begin{equation} \label{pulch:truncation}
  \tilde{y}^{({I})}(t,{\bf p}) =
  \sum_{i \in {I}} \tilde{w}_i(t) \Phi_i({\bf p})
\end{equation}
with a finite index set ${I} \subset \mathbbm{N}$
and approximations~$\tilde{w}_i$ of the coefficients.
Typically, all polynomials up to a total degree~$d$ are included
in an index set ${I}^d$.
The number of basis polynomials becomes
$\left| {I}^d \right| = \frac{(d+q)!}{d!q!}$.

Stochastic collocation techniques yield approximations of the
unknown coefficient functions, see~\cite{pulch:xiu-book}.
A quadrature rule is determined by
nodes $\{ {\bf p}^{(1)} , \ldots , {\bf p}^{(s)} \} \subset \Pi$
and weights $\{ \gamma_1,\ldots, \gamma_s \} \subset \mathbbm{R}$.
The approximations become
\begin{equation} \label{pulch:collocation}
\tilde{w}_i(t) = 
\sum_{j=1}^s \gamma_j \, y(t,{\bf p}^{(j)}) \, \Phi_i({\bf p}^{(j)})
\end{equation}
for $i=1,\ldots,m$ and w.l.o.g. ${I}^d = \{ 1,\ldots,m \}$.
Thus the coupled problem~(\ref{pulch:mna}),(\ref{pulch:pde})
has to be solved $s$-times for different realisations of the parameters. 


\section{Sparse Approximation}
\label{pulch:sec:sparse}
The aim is to find an index set ${J} \subset {I}^d$
with $| {J} | \ll | {I}^d |$ for fixed total degree~$d$,
while the error is still below some threshold.
The total error $y - \tilde{y}^{({J})}$ consists of three parts:
(i) the truncation error ($\mathbbm{N} \rightarrow {I}^d$),
(ii) the error of the numerical method ($w_i \rightarrow \tilde{w}_i$), and
(iii) the additional sparsification error
(${I}^d \rightarrow {J}$).
We assume that the errors (i) and (ii) are sufficiently small and
focus on the error~(iii).

The relative ${L}^2$-error of the sparsification reads as
\begin{equation} \label{pulch:error}
  E(t;{J}) =
  \frac{\left\| \tilde{y}^{({I}^d)} (t, \cdot ) 
  - \tilde{y}^{({J})} (t, \cdot ) \right\|_{{L}^2}}{
  \left\| \tilde{y}^{({I}^d)} (t, \cdot ) \right\|_{{L}^2}} =
\left( \frac{  \displaystyle
  \sum_{i \in {I}^d \backslash {J}} \tilde{w}_i(t)^2 }{
  \displaystyle \sum_{i \in {I}^d} \; \tilde{w}_i(t)^2 } \right)^{\frac{1}{2}}
\end{equation}
for each~$t$ including the coefficients~(\ref{pulch:collocation}).
Given an error tolerance $\varepsilon > 0$,
we obtain an optimal index set 
\begin{equation} \label{pulch:optimal-index}
  {J}_t = {\rm argmin} 
  \left\{ \left| {J}' \right| \; : \;
  {J}' \subseteq {I}^d  \;\; \mbox{and} \;\;
  E(t,{J}') < \varepsilon \right\}
\end{equation}
with respect to the error~(\ref{pulch:error}) for each time point.
A global index set is given by
\begin{equation} \label{pulch:global-index}
  \hat{J} = \displaystyle \bigcup_{t \in [t_0,t_{\rm end}]} {J}_t ,
\end{equation}
which is sufficiently accurate with respect to the tolerance~$\varepsilon$
for all times.
More details can be found in~\cite{pulch:matcom2018}.


\section{Model Order Reduction}
\label{pulch:sec:pod}
Alternatively, we use an MOR with proper orthogonal decomposition (POD),
see~\cite{pulch:antoulas}, to determine a low-dimensional approximation
of the polynomial surrogate.
Let ${I}^d = \{ 1,\ldots,m \}$ and
$\tilde{\bf w} = (\tilde{w}_1,\ldots,\tilde{w}_m)^\top$.
A transient simulation of the coupled problem yields
the coefficients~(\ref{pulch:collocation}) by the
stochastic collocation technique.
We collect snapshots
$\tilde{\bf w}(t_0),\tilde{\bf w}(t_1),\ldots,\tilde{\bf w}(t_{k-1})$
for discrete time points 
in a matrix ${\bf W} \in \mathbbm{R}^{m \times k}$.
A singular value decomposition yields
$$ {\bf W} = {\bf U} {\bf S} {\bf V}^\top
\qquad \mbox{with} \qquad
{\bf S} = {\rm diag}(\sigma_1,\sigma_2,\ldots,\sigma_{\min\{m,k\}}) $$
including the singular values $\sigma_1 \ge \sigma_2 \ge \cdots$
and orthogonal matrices ${\bf U} \in \mathbbm{R}^{m \times m}$,
${\bf V} \in \mathbbm{R}^{k \times k}$.
Let ${\bf u}_1,\ldots,{\bf u}_m$ be the columns of the matrix~${\bf U}$.
We arrange the orthogonal projection matrix 
${\bf P}_r = \left( {\bf u}_1 \, \cdots \, {\bf u}_r \right)
\in \mathbbm{R}^{m \times r}$
(${\bf P}_r^\top {\bf P}_r = {\bf I}_r$)
for each dimension $r \le \min \{ m,k \}$.
Given $\tilde{\bf w}(t) \in \mathbbm{R}^m$,
the best approximation with respect to the reduced basis reads as 
$\bar{\bf w}_r(t) = {\bf P}_r^\top \tilde{\bf w}(t)$
for any~$t$.
Vice versa, we obtain the approximation
$\tilde{\bf w}(t) \approx {\bf P}_r \bar{\bf w}_r(t)$
for given $\bar{\bf w}_r(t)$ and any~$t$.

We define the (dense) low-dimensional approximation,
cf.~(\ref{pulch:truncation}),
\begin{equation} \label{pulch:low-dimensional}
  \tilde{y}^{({I}^d)}(t,{\bf p}) \approx \displaystyle
  \sum_{i = 1}^m \left[ \sum_{j=1}^r u_{ij} \bar{w}_j(t) \right] \Phi_i({\bf p})
  = \sum_{j=1}^r \bar{w}_j(t)
  \underbrace{\left[ \sum_{i=1}^m u_{ij} \Phi_i({\bf p}) \right]}_{=:\Psi_j({\bf p})} 
\end{equation}
with the new orthonormal basis polynomials $\{ \Psi_1 , \ldots , \Psi_r \}$
and associated coefficients $\bar{w}_1,\ldots,\bar{w}_r$.
An error estimate for an approximation of this kind is given
in~\cite{pulch:arxiv}.


\section{Numerical Results}
\label{sec:results}
In the coupled problem~(\ref{pulch:mna}),(\ref{pulch:pde}),
we supply a harmonic oscillation with period $T=0.02$ as input voltage.
We use the Stroud-5 quadrature rule with $s=243$ nodes,
which is exact for all polynomials up to total degree~5,
see~\cite{pulch:stroud}.
In each node, we perform a monolithic time integration in $[0,2T]$
by the implicit Euler method.
The step size $\Delta t = 10^{-4}$ is used in time,
whereas a smaller step size does hardly change the numerical results.
This time integration yields $k=401$ snapshots in equidistant time points.
We choose the polynomial degree $d=3$, i.e.,
$\left| {I}^3 \right| = m = 364$ due to $q=11$ random parameters.

\begin{figure}
  \begin{center}
    \includegraphics[width=5.5cm]{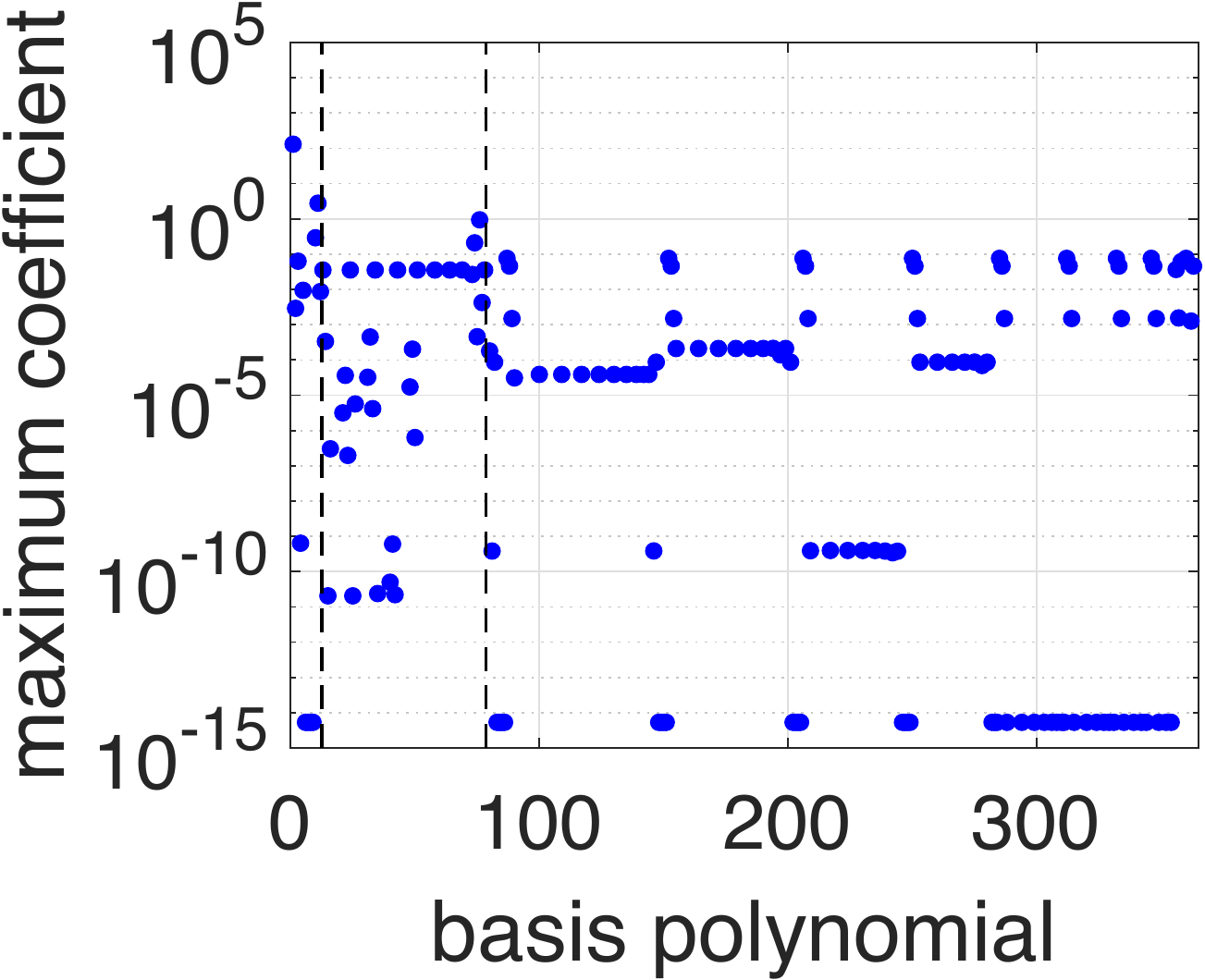}
    \hspace{5mm}
    \includegraphics[width=5.5cm]{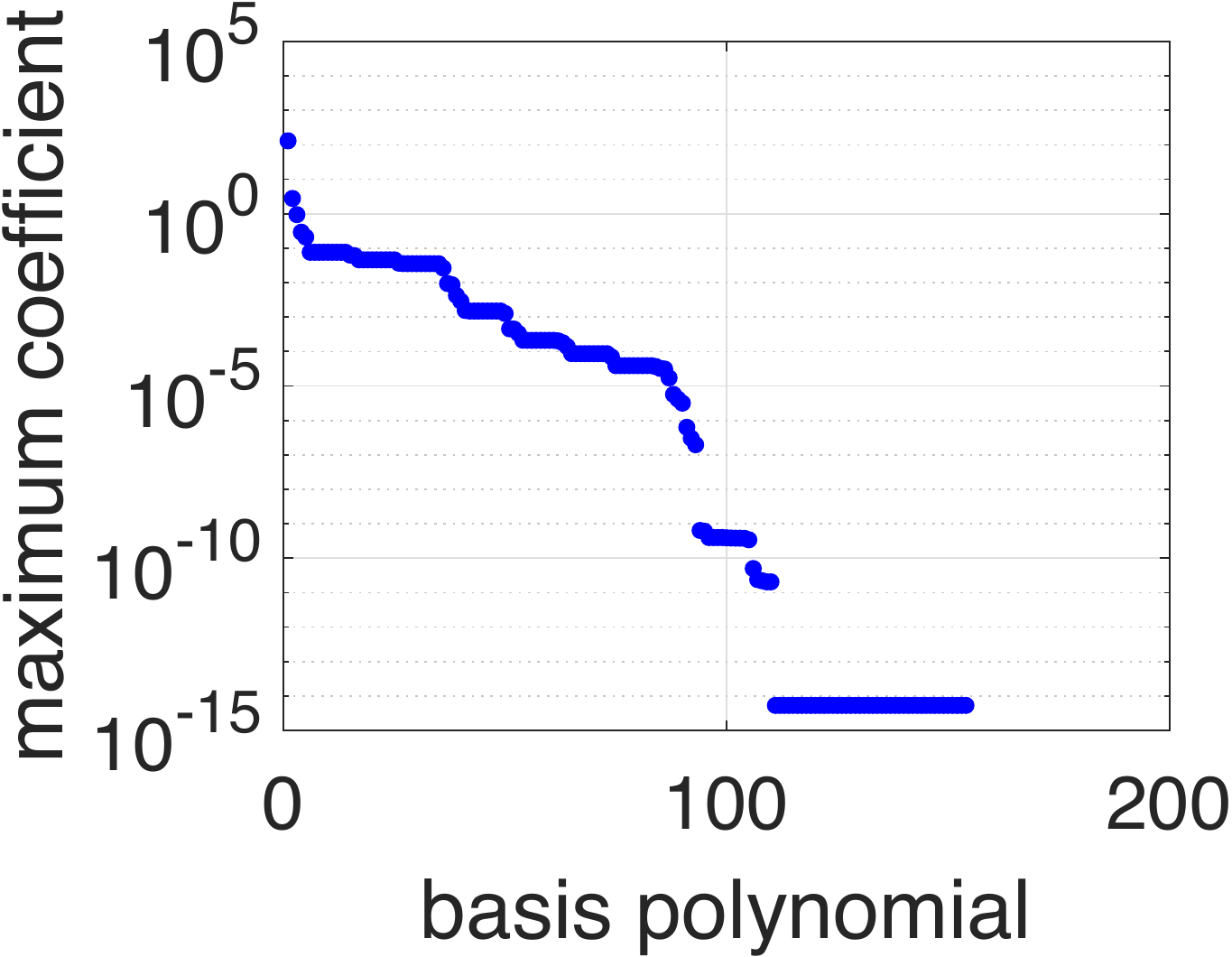}
  \end{center}
  \caption{Maximum of $\{ \tilde{w}_i(t) \; : \; t \in [0,2T] \}$
    for $i=1,\ldots,364$ with coefficients~(\ref{pulch:collocation}),
    left: dashed lines separate the coefficients of degree zero/one,
    two and three,
    right: coefficients in descending order.}
  \label{pulch:fig:coefficients}
\end{figure}

\begin{figure}
  \begin{center}
    \includegraphics[width=5.5cm]{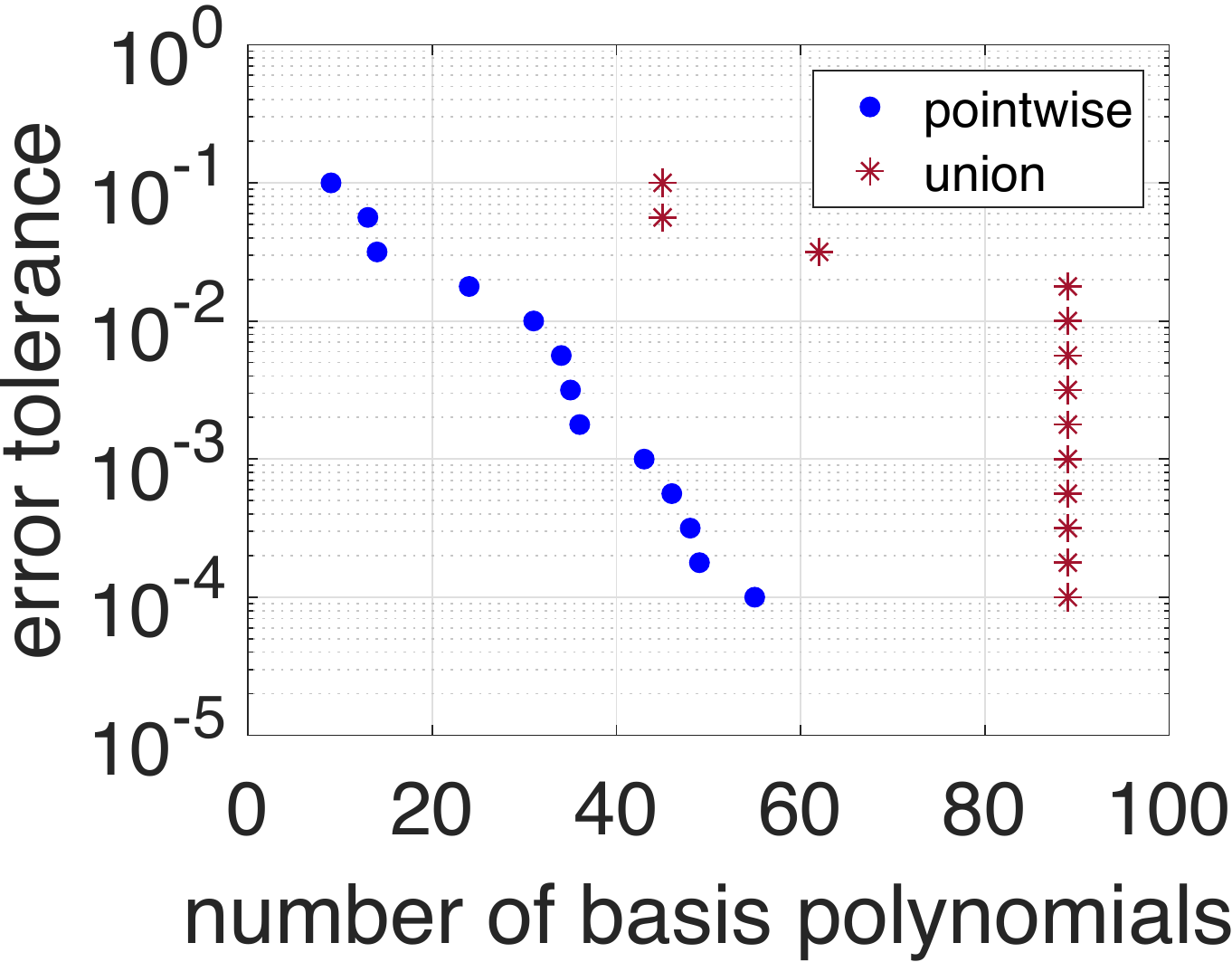}
    \hspace{5mm}
    \includegraphics[width=5.5cm]{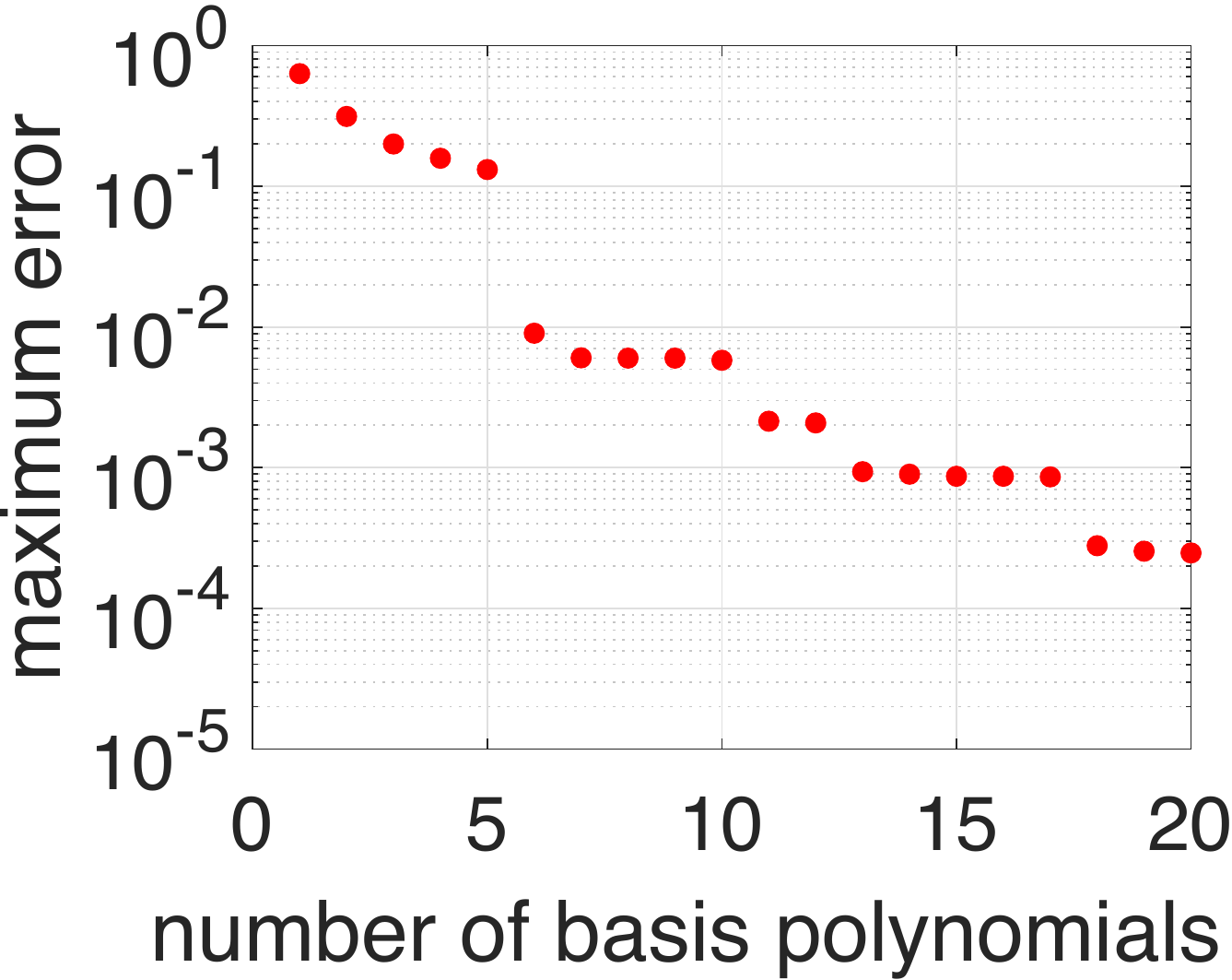}
  \end{center}
  \caption{Relation between number of basis polynomials and error tolerance
    or maximum relative ${L}^2$-errors,
    left: sparsification by neglecting basis polynomials,
    right: basis selection using POD.}
  \label{pulch:fig:error}
\end{figure}

Figure~\ref{pulch:fig:coefficients} (left) shows the maximum
of the coefficients~(\ref{pulch:collocation}) in time.
All coefficients of degree three are (at least) three orders of magnitudes
smaller than the coefficient of degree zero.
This property suggests that the relative truncation error is below
0.1\%.
Furthermore, Figure~\ref{pulch:fig:coefficients} (right)
demonstrates a fast decay,
which indicates some potential for a sparse approximation
as described in Section~\ref{pulch:sec:sparse}.
For given error tolerances $\varepsilon \in [10^{-4},10^{-1}]$,
we determine the cardinalities
$\max_{t \in [0,2T]} | {J}_t |$ (pointwise)
and $| \hat{{J}} |$ (union) with the index
sets~(\ref{pulch:optimal-index}) and~(\ref{pulch:global-index}),
respectively.
Figure~\ref{pulch:fig:error} (left) illustrates the
cardinalities in dependence on the error tolerances.
The potential for a global sparse approximation
using~(\ref{pulch:global-index}) is bad,
because more than 80~basis polynomials are required.

Alternatively, we apply the POD technique from Section~\ref{pulch:sec:pod}.
We choose the reduced dimensions $r=1,\ldots,20$ and obtain
the approximations~(\ref{pulch:low-dimensional}).
Figure~\ref{pulch:fig:error} (right) depicts the maximum in time
of the relative ${L}^2$-errors for these approximations.
Now we achieve an efficient low-dimensional representation,
where less than 20~basis polynomials ($r \ll m$) yield a small error.

\section{Conclusions}

We performed an UQ of a coupled DAE-PDE system modelling
a field-circuit problem.
Sparse approximations of the random QoI were identified
using its orthogonal polynomial expansion.
In this test example, an appropriate sparse representation could not
be achieved on the global time interval by simply
neglecting basis polynomials.
Alternatively, we obtained an efficient low-dimensional approximation
by changing to another orthogonal polynomial basis,
which was identified via an MOR.

\begin{acknowledgement}
  The research of the second author is supported by the
  ‘Excellence Initiative’ of the German Federal and State Governments and
  by the Graduate School of Computational Engineering at Technische
  Universit\"at Darmstadt.
\end{acknowledgement}


\end{document}